\documentstyle[12pt]{amsart}

\makeatletter                                            
\renewcommand{\@begintheorem}[2]{                        
\rm \trivlist \item [\hskip \labelsep {\bf #2\ \ #1.}]   
                                }                        
\makeatother                                             

\newcommand{\newsubsection}%
{{\bf\refstepcounter{subsection}\thesubsection\ \ }}

\newcommand{\newsubsubsection}%
{{\bf\refstepcounter{subsubsection}\thesubsubsection\ \ }}

\newcommand{\ts}{\vspace{\baselineskip}{\bf Proof.}$\;\;$}
\newcommand{\ZZ}{{\bf Z}}
\newcommand{\QQ}{{\bf Q}}
\newcommand{\RR}{{\bf R}}

\newcommand{\CC}{{\bf C}}
\newcommand{\FF}{{\bf F}}
\newcommand{\HH}{{\bf H}}
\newcommand{\PP}{{\bf P}}

\newcommand{\cA}{{\cal A}}

\newcommand{\cI}{{\cal I}}

\newcommand{\cM}{{\cal M}}

\newcommand{\cO}{{\cal O}}

\newcommand{\hC}{{\hat{C}}}
\newcommand{\tC}{{\tilde{C}}}
\newcommand{\tX}{{\tilde{X}}}
\newcommand{\tx}{{\tilde{x}}}

\begin{document}

\title{Quaternionic pryms and Hodge classes}
\author{B. van Geemen}
\address{Dipartimento di Matematica, Universit\`a di Pavia,
  via Ferrata 1, I-27100 Pavia, Italia}
\email{geemen@@dragon.ian.pv.cnr.it}
\author{A. Verra}
\address{Dipartimento di Matematica, Universit\`a di Roma 3, 
Largo San Leonardo Murialdo 1
I-00146 Roma, Italia}
\email{verra@@mat.uniroma3.it}

\begin{abstract}
Abelian varieties of dimension $2n$ on which a definite quaternion algebra acts
are parametrized by symmetrical domains of dimension $n(n-1)/2$. Such abelian 
varieties have primitive Hodge classes in the middle dimensional cohomology group. In general, it is not clear that these are cycle classes. In this paper we show that a particular $6$-dimensional family of such $8$-folds are Prym varieties and we use the method of C.\ Schoen to show that all Hodge classes on the general abelian variety in this family are algebraic. 
We also consider Hodge classes on certain $5$-dimensional
subfamilies and relate these to the Hodge conjecture for abelian $4$-folds.
\end{abstract}

\maketitle

In this paper we study abelian varieties of dimension $8$ whose endomorphism ring is a definite quaternion algebra over $\QQ$, we refer to these as abelian $8$-folds 
of quaternion type. Such abelian varieties are interesting 
since their Hodge rings are not generated by divisor classes \cite{KM} 
and the Hodge conjecture is still open for most of them. 
The moduli space of $8$-folds of quaternion type, with fixed discrete data,
is $6$-dimensional. 
One such moduli space was investigated recently by Freitag and Hermann
\cite{FrHe}. In \cite{shiga} a 6-dimensional family of K3 surfaces is studied
whose Kuga-Satake varieties have simple factors of dimension $8$ which are of quaternion type. 

Our first result, in section \ref{prym surj}, is that 
one of these moduli spaces
parametrizes  Prym varieties of unramified 2:1 covers 
$\tC\rightarrow \hC$, actually the genus $17$ curves $\tC$ are 8:1 unramified 
covers of genus $3$ curves with Galois group $Q$ generated by the quaternions
$i$ and $j$.  
Results of C.\ Schoen \cite{S1}
on cycles on Prym varieties imply that the Hodge conjecture holds 
for the general quaternionic $8$-fold in this family, see Corollary \ref{cycler}.
Previously, S. Abdulali \cite{A} showed the Hodge conjecture for certain 
$3$-dimensional families of quaternionic $6$-folds using results on Pryms of 
cyclic 3:1 covers. 

The remainder of the paper is devoted to the relation between the
Hodge cycles on a quaternionic $8$-fold and on a self-product of an abelian $4$-fold of Weil type. This is inspired by Schoen's paper \cite{S2} and by
the fact that the Hodge conjecture for all abelian $4$-folds follows
 {\em if} for all abelian $4$-folds of 
Weil type the Weil classes are classes of algebraic cycles \cite{MZ}. We recall the definition of these varieties and classes in section \ref{weil}.

The self-product of an abelian variety of Weil type is of quaternion type.
Our Proposition \ref{so8res} shows that the exceptional Hodge classes on a general quaternionic variety specialize to classes of intersections of divisors on a self-product.
This is unfortunate, since till now nothing is known about he Hodge conjecture for abelian $4$-folds of Weil type with fields $K$ which are distinct from
$\QQ(\sqrt{-1})$ and $\QQ(\sqrt{-3})$. 

In the last section we study certain $5$-dimensional families of $8$-folds of quaternion type, called abelian varieties of $so(7)$-type. 
We show that the Hodge $(2,2)$-conjecture for the general member of such a family implies the Hodge conjecture for infinitely many ($4$-dimensional) families of abelian varieties of Weil type with field $K$ where $K$ can be any subfield of the quaternion algebra.
We hope to identify such $5$-dimensional families explicitly (in terms of the 
Prym construction) and to investigate the Hodge $(2,2)$ conjecture for the corresponding Pryms in the future.  

We are indebted to G.\ Lombardo for some interesting discussions.

\section{Quaternionic covers}

\subsection{Definitions.}\label{defs}
A definite quaternion algebra $F$ over $\QQ$ is a skew field with center $\QQ$,
of dimension 4 over $\QQ$, such that $F\otimes_\QQ\RR$ is isomorphic to the skew field of Hamilton's quaternions. It is well known that one can find elements
$i,\,j,\,k\in F$ such that
$$
F=\QQ+\QQ i+\QQ j+\QQ k,\qquad i^2=r,\quad j^2=s,\quad ij=-ji=k,
\qquad{\rm with}\quad r,\;s\in\QQ_{<0}.
$$
Such a quaternion algebra has a canonical (anti-)involution, it acts as: 
$$
x=a+bi+cj+dk\longmapsto \bar{x}=a-bi-cj-dk,\qquad \overline{xy}=\bar{y}\bar{x}.
$$
The quaternion algebra over $\QQ$ with $r=s=-1$ will be denoted by $\HH_\QQ$.
The quaternion group $Q$ is the following subgroup of order 8 of $\HH_\QQ^*$:
$$
Q=\{\pm 1,\pm i,\;\pm j,\;\pm k\},\qquad (\pm i)^2=(\pm j)^2=(\pm k)^2=-1,\quad
ij=-ji=k.
$$

\subsection{Quaternionic covers.}
Let $\cM_{g,Q}$ be the moduli space of unramified Galois covers $\tilde{C}\rightarrow C$ of a genus $g$ curves $C$ with Galois group $Q$. 
A quaternionic cover $\pi:\tilde{C}\rightarrow C$ defines a normal subgroup $N:=\pi_*\pi_1(\tC)$ of $\pi_1(C)$ and $\pi_1(C)/N\cong Q$.
Conversely, a normal subgroup $N$ of $\pi_1(C)$ with $\pi_1(C)/N\cong Q$
defines a quaternionic cover:
$$
N\backslash C^{univ}\longrightarrow C,
$$
where $C^{univ}$ is the universal cover of $C$.
Since there exist surjective homomorphisms $\pi_1(C)\rightarrow Q$, 
there is a surjective map $\cM_{g,Q}\rightarrow {\cM}_g$ where 
$\cM_g$ is the moduli space of genus $g$ curves. This map is finite since there
the number of surjective homomorphisms $\pi_1(C)\rightarrow Q$ is finite 
(in fact, $\pi_1(C)$ is finitely generated).

Proposition \ref{qstan} shows that, up to automorphisms of $\pi_1(C)$, there is only one subgroup defining a quaternionic cover. Standard monodromy arguments then imply that the moduli space $\cM_{g,Q}$ is irreducible.

\subsection{}
The fundamental group of a genus $g$ curve is the quotient of the free group $F_{2g}$ generated by $a_1,\ldots,a_g$, $b_1,\ldots,b_g$ by the normal subgroup generated by the relation $R=[a_1,b_1]\ldots[a_g,b_g]$. Let $A_g$ be the subgroup of $Aut(F_{2g})$ of automorphisms which map $R$ to a conjugate of itself, there is a natural map
$A_g\rightarrow Aut(\pi_1(C))$.
The mapping class group $\Gamma_g$ is the quotient of $A_g$ by inner automorphisms
of $\pi_1(C)$.

\subsection{Proposition.}\label{qstan}
Let $N\subset\pi_1(C)$ be a normal subgroup such that $\pi_1(C)/N\cong Q$.
Then there exists an automorphism $\psi\in A_g$ such that
$$
N=\ker(f\circ\psi :\pi_1(C)\longrightarrow Q)
$$
with $f$ the homomorphism:
$$
f:\pi_1(C)\longrightarrow Q\qquad f(\alpha_{g-1})=i,\quad f(\alpha_g)=j,\qquad f(\alpha_i)=f(\beta_j)=1
$$
for all $i$ and $j$ with $1\leq i\leq g-2$ and $1\leq j\leq g$
where $\alpha_j$, $\beta_j$ are generators of $\pi_1(C)$ satisfying
$[\alpha_1,\beta_1]\ldots [\alpha_g,\beta_g]=e$.

\ts
The composition 
$$
\pi_1(C)\longrightarrow \pi_1(C)/N\cong Q\longrightarrow Q/\{\pm 1\}\cong
(\ZZ/2\ZZ)^2
$$
factors over $H_1(C,\ZZ/2\ZZ)$ and the group $A_g$ maps onto the 
symplectic group (w.r.t.\ the intersection form) of this $\FF_2$-vector space. 
Therefore there is a $\psi\in A_g$ such that $N=\ker(g\circ \psi)$
with a homomorphism $g:\pi_1(C)\rightarrow Q$ given by $g(\alpha_g)=j$
and either $g(\alpha_{g-1})=i$ (and $g$ maps the other generators to $\{\pm 1\}$)
or $g(\beta_g)=i$ (and $g$ maps the other generators to $\{\pm 1\}$).
However, if $g(\beta_g)=i$ then $[\alpha_g,\beta_g]=-1$ whereas 
$[\alpha_i,\beta_i]=1$ for $i<g$, which contradicts the relation $R$. 
Therefore $g(\alpha_{g-1})=i$. 

It remains to show that, using elements of $A_g$, we may assume that $g$ maps 
the other generators to $+1$. If $g(\beta_g)=-1$, 
then we apply the automorphism $\beta_g\mapsto \beta_g\alpha_g^2$ 
(and fixing the other generators), similarly if $g(\beta_{g-1})=-1$. 

If for $i<g-1$ we have $g(\alpha_i)=-1$ or $g(\beta_i)=-1$, 
we use the automorphism in $A_g$ which fixes all generators except $\alpha_i$ 
which is mapped to
$\alpha_i\mapsto \alpha_i\beta_i$ and/or the automorphism in $A_g$ 
which fixes all generators except $\beta_i$ which is mapped to 
$\beta_i\alpha_i$ to obtain
a homomorphism  which maps $\alpha_i\mapsto -1$, $\beta_i\mapsto 1$.

Thus we are reduced to the case of a homomorphism $f:\pi_1(C)\rightarrow Q$
which is trivial on the $\beta_i$, maps the $\alpha_i$ to $\pm 1$ for $i<g-1$
and has $f(\alpha_{g-1})=i$ and $f(\alpha_g)=j$.
For an integer $k$ we define $\psi_k$ as follows:
$$
\psi_k:\left\{
\begin{array}{rcl}
\alpha_i&\longmapsto&\alpha_{k+1}\alpha_i\alpha_{k+1}^{-1}\\
\beta_i&\longmapsto&\alpha_{k+1}\beta_i\alpha_{k+1}^{-1}\\
\end{array}\qquad
\begin{array}{rcl}
\alpha_k&\longmapsto&\alpha_{k+1}\alpha_k\\
\beta_k&\longmapsto&\beta_k\\
\end{array}\qquad
\begin{array}{rcl}
\alpha_{k+1}&\longmapsto&\beta_k\alpha_{k+1}\beta_k^{-1}\\
\beta_{k+1}&\longmapsto&\alpha_{k+1}\beta_{k+1}\alpha_{k+1}^{-1}\beta_k^{-1}\\
\end{array}\right.
$$
where $1\leq i\leq g$ and $i\neq k,\,k+1$.
It is easy to see that $\psi_k$ is an automorphism of the free group generated 
by the $\alpha_j$ and $\beta_j$, and since 
$
\psi_k([\alpha_i,\beta_i])=\alpha_{k+1}[\alpha_i,\beta_i]\alpha_{k+1}^{-1}$
and
$$
\psi_k([\alpha_k,\beta_k])=
\alpha_{k+1}[\alpha_k,\beta_k]\beta_k\alpha_{k+1}^{-1}\beta_k^{-1},\qquad
\psi_k([\alpha_{k+1},\beta_{k+1}])=
\beta_k\alpha_{k+1}\beta_k^{-1}[\alpha_{k+1},\beta_{k+1}]\alpha_{k+1}^{-1}
$$
it follows that $\psi_k(R)=\alpha_{k+1}R\alpha_{k+1}^{-1}$, 
hence $\psi_k\in A_g$. 
Note that $(f\circ\psi_k)(\beta_i)=1$ for all $i$ and $k$. 
Moreover, $(f\circ \psi_{g-2}^2)(\alpha_{g-2})=-f(\alpha_{g-2})$ and 
$f\circ \psi_{g-2}^2$ has the same value as $f$ on the other generators. 
Similarly, for $k<g-2$, we have
$(f\circ \psi_{k})(\alpha_{k})=f(\alpha_{k})f(\alpha_{k+1})$ $=\pm f(\alpha_k)$
and $f\circ \psi_{k}$ has the same value on the other generators. 
Therefore composing $f$ with a suitable composition of $\psi_k$'s we get the 
homomorphism as desired in the proposition.
\qed

\subsection{Corollary.} \label{qirr}
The moduli space $\cM_{g,Q}$ of quaternionic covers is irreducible.

\ts
Recall that the Teichm\"uller space ${\cal T}_g$ parametrizes genus $g$ curves with a set of standard generators of $\pi_1(C)$. The group $\Gamma_g$ acts on 
${\cal T}_g$ and $\cM_g=\Gamma_g\backslash{\cal T}_g$.
Proposition \ref{qstan} shows that the map
$$
{\cal T}_g\longrightarrow \cM_{g,Q},\qquad [C,\{\alpha_1,\ldots,\beta_g\}]
\longmapsto [(\ker{f})\backslash C^{univ}\rightarrow C],
$$
where $f$ is as in the proposition, is surjective, hence $\cM_{g,Q}$ is irreducible.
\qed

\subsection{Construction of quaternionic covers.}
For completeness sake we briefly indicate how one can construct quaternionic covers of an algebraic curve $C$ of genus $g>1$. By Corollary \ref{qirr}, any quaternionic cover arises in this way.

Choose an \'etale double cover $\pi_\alpha:C_\alpha\rightarrow C$ of degree 2 with covering involution $\iota_\alpha$. 
Let $\eta\in Div(C_\alpha)$ be a divisor such that $4\eta\equiv 0$ and $\eta+\iota_\alpha^*\eta\equiv 0$ and let $h\in k(C_\alpha)$, the function field of $C_\alpha$, have divisor $(h)=\eta+\iota_\alpha^*\eta$.
Then $\iota_\alpha^*h=\pm h$ and we claim that, changing maybe $\eta$, we may assume that $\iota_\alpha^*h=-h$.

To prove the claim we choose $\gamma\in Div(C)$ such that $2\gamma\equiv \alpha$, a divisor of order 2 in $Pic(C)$ which defines $C_\alpha$.
Then $\pi_\alpha^*\gamma$ has order two in $Pic(C_\alpha)$ and $\iota_\alpha^*\pi_\alpha^*\gamma=\pi^*_\alpha\gamma$ whence $\iota_\alpha^*\pi_\alpha^*\gamma+\pi^*_\alpha\gamma\equiv 0$.
Therefore there is a rational function $g\in k(C_\alpha)$ with $(g)=
\iota_\alpha^*\pi_\alpha^*\gamma+\pi^*_\alpha\gamma$ and $\iota_\alpha^*g=-g$ since the divisor of $g$ is invariant under $\iota^*$ but $g\not\in \pi_\alpha^*k(C)\;(\subset k(C_\alpha))$. Now in case $\iota_\alpha^*h=h$, we let
$\eta':=\eta+\pi_\alpha^*\gamma$. Then $4\eta'\equiv 0$,
$\eta'+\iota_\alpha^*\eta'=(gh)\equiv 0$ and moreover $\iota_\alpha^*(gh)=-gh$ which proves the claim.

Now consider the cover $\tilde{C}$ of $C_\alpha$ which is defined by the field
$k(\tilde{C}):=k(C_\alpha)[T]/(T^4-f)$ where $f\in k(C_\alpha)$ is such that
$(f)=4\eta$. Thus $\tilde{C}$ is an \'etale 4:1 cover of $C_\alpha$.
We claim that the \'etale cover $\tilde{C}\longrightarrow C$ is quaternionic.

To see this, let $t:=T+(T^4-f)\in k(\tilde{C})$ and define 
$$
\tilde{\iota}^*:k(\tilde{C})\longrightarrow k(\tilde{C}),\qquad
\tilde{\iota}^*(t)=\frac{h}{t},\quad \tilde{\iota}^*_{|k(C_\alpha)}=\iota_\alpha^*.
$$
Then $\tilde{\iota}^*$ is an automorphism of $k(\tilde{C})$ which lifts $\iota^*\alpha$
since $f\iota_\alpha^*f=h^4$ implies
$$
\tilde{\iota}^*(t^4-f)=\tilde{\iota}^*(t)^4-\iota_\alpha^*f=
(h/t)^4-\iota_\alpha^*f=h^4/f-\iota_\alpha^*f=0.
$$
The order of the automorphism $\tilde{\iota}^*$ is 4 since 
$$
\tilde{\iota}^*(\tilde{\iota}^*(t))=\tilde{\iota}^*(h/t)=\tilde{\iota}^*(h)t/h=
-t.
$$
Any other lift of $\iota_\alpha^*$ is given by $t\mapsto i^ah/t$ for $a=1,\,2,\,3$ (and $i^2=-1$) hence all lifts of $\iota_\alpha^*$ to $k(\tilde{C})$ have order 4
and therefore the cover $\tilde{C}\rightarrow C$ is quaternionic.

\section{Quaternionic pryms.}

\subsection{Quaternion type.}
We will say that $(A,E,F)$ is a polarized abelian variety of quaternion type if
$(A,E)$ is a polarized abelian variety and if 
there is a definite quaternion algebra $F$ over $\QQ$ with 
$F\subset End(A)\otimes \QQ$ in such a way that $H_1(A,\QQ)$ is a vector space over $F$ and $x^*E=x\bar{x} E$. 
The last condition is equivalent to the condition
that the Rosati involution defined by $E$ on $End(A)$ induces the canonical involution on $F$.

\subsection{Moduli.}\label{moduli}
Let now $(A,E,F)$ be a polarized abelian variety of quaternion type, 
and let $\dim A=2n$. The deformations of $A$ which have the same endomorphism ring and polarization
are parametrized by a symmetric space which has dimension $n(n-1)/2$ (cf.\ \cite{LB}, section 9.5).

\subsection{Quaternionic covers.}\label{qp}
Let $\tC\rightarrow C$ be a quaternionic cover of a genus $g$ curve $C$.
 The Zeuthen-Hurwitz formula shows that the genus of $\tilde{C}$ is $8g-7$. The quotient of $\tilde{C}$ by $\{\pm 1\}\subset Q$ is a curve of genus $4g-3$, denoted by $\hat{C}$.
$$
\tilde{C}\stackrel{2:1}{\longrightarrow} \hat{C}=\tilde{C}/\{\pm 1\} \stackrel{4:1}{\longrightarrow} C.
$$
The Prym variety $P:=Prym(\tilde{C}/\hat{C})$ is called a quaternionic Prym.
The following proposition shows that it is of quaternion type.

\subsection{Proposition.} 
The Prym variety $P:=Prym(\tilde{C}/\hat{C})$ is a $4(g-1)$-dimensional 
principally polarized abelian variety of quaternion type for the algebra 
$\HH_\QQ$.

\ts
The dimension of $P$ is $8g-7-(4g-3)=4(g-1)$ and Pryms of unramified double covers are principally polarized. The group $Q$ has 5 irreducible complex representations, 4 of which are one-dimensional, these factor over $Q/\{\pm 1\}\cong (\ZZ/2\ZZ)^2$, and one is faithful and has dimension 2.
The group ring $\QQ[Q]$ is isomorphic to $\QQ^4\oplus \HH_\QQ$.
By definition of the Prym, $-1\in Q$ acts as $-1$ on the tangent space $T_0P$
hence $H_1(P,\QQ)$ is a $\HH_\QQ$-vector space. 
Let $\alpha\in\{i,\,j,\, k\}\subset Q$,
then $\alpha^*\in Aut(P)$ and thus $E(\alpha^*x,\alpha^*y)=E(x,y)$ for all
$x,\,y\in H_1(P,\QQ)$ where $E\in \wedge^2H^1(P,\QQ)$ is the polarization on $P$. Since the Rosati involution $\phi\mapsto \phi'$ is characterized by $
E(\phi x,y)=E(x,\phi'y)$ and $\alpha^2=-1$ we conclude that $\alpha'=-\alpha^*$.
As $(a\phi+b\psi)'=a\phi' +b\psi'$ for $a,\,b\in\QQ$ and $\phi,\,\psi\in End(P)$, the Rosati involution induces the canonical involution on $\HH_\QQ$.
\qed

\subsection{} We determine the homomorphism
$$
Q\longrightarrow Aut(H_1(P,\ZZ))\cong Sp(4(g-1),\ZZ).
$$
The period matrices which correspond to 
quaternionic Pryms are contained in the fixed point set of the action of $Q$ 
(via the symplectic group) on the Siegel upper half space.
To do this we introduce the ring
$$
\HH_\ZZ:=\{a+bi+cj+dk\in\HH_\QQ:\;a,\,b,\,c,\,d\in\ZZ\;\}.
$$ 
and the $\HH_\ZZ$-module (actually itself a ring, 
the ring of Hurwitz integers):
$$
M:=\{a\zeta+bi+cj+dk\in\HH_\QQ:\;a,\,b,\,c,\,d\in\ZZ\;\},\qquad
{\rm with}\quad \zeta=(1+i+j+k)/2.
$$
The $\HH_\ZZ$-module $M$ is not free since there is no $m\in M$ such that $m,\;im,\;jm,\;km$ are a $\ZZ$-basis of $M$, in fact, $(1+i+j+k)m$
is divisible by $2$  in $M$ for any $m\in M$.

\subsection{Proposition.}
The period lattice of a quaternionic Prym of dimension $4(g-1)$ 
is isomorphic to the $\HH_\ZZ$-module:
$$
\Lambda\cong \left(M\oplus \HH_\ZZ^{g-2}\right)^2.
$$

\ts
We consider first the case that $C$ has genus two, we may assume that the cover $\tC\rightarrow C$ is defined by $\ker(f)$ with $f$ as in Proposition \ref{qstan}. Let $X$ be the union of two circles with a point $x_0$ in common,
note that $\pi_1(X)$ is a free group on two generators. 
We can view $C$ as the boundary of a tubular neighbourhood of $X$ and we may assume that the contraction $C\rightarrow X$ maps $\alpha_1$, $\alpha_2\in \pi_1(C)$ homeomorphically onto the two circles. Then $\pi:\tC\rightarrow C$ is
the pull-back of a cover $p:\tX\rightarrow X$. 
Choosing a base point $\tx_0\in p^{-1}(x_0)$, we can identify 
$Q\simeq p^{-1}(x_0)$ via $g\mapsto  [g]:=g\cdot\tx_0$. 
As $f(\alpha_1)=i$, $f(\alpha_2)=j$, 
the 8 distinct lifts of the first circle are paths in $\tX$ connecting the points $g$ and $ig$, whereas the lifts of the second circle
connect $g$ and $jg$. 
In particular, $\tX$ is a simplicial complex with vertices $[g]$, $g\in Q$, 
and edges $[g,ig]$ and $[g,jg]$. 
Let 
$$
c=[1,i]+[i,-k]+[-k,-j]+[-j,1]+[-i,-1]+[-1,j]+[j,k]+[k,-i].
$$
Then one easily verifies that the subgroup $V_-$ of $H_1(\tX,\ZZ)$ 
generated by the cycles $c$, $ic$, $jc$ and $\zeta c$ is a direct summand of 
$H_1(\tX,\ZZ)$. In particular, $V_-\cong M$. 
Since $-1\in Q$ acts as $-1$ on $V_-$, we may identify $M$ with a direct summand
of $H_1(P,\ZZ)\subset H_1(\tC,\ZZ)$ 
via the contraction $\tC\rightarrow \tX$.
Note that the intersection pairing is trivial on $V_-\times V_-$.

The inverse image of an open edge $e^0\subset\tX$ is a cylinder 
$S^1\times ]0,1[$
in $\tC$ and we define $b_e\in H_1(\tC,\ZZ)$ to be the homology class of the 
loop $S^1\times \{1/2\}$. 
The subgroup $H_1(P,\ZZ)$ of $H_1(\tC,\ZZ)$ on which $-1\in Q$ acts as $-1$
is the direct sum of $V_-$ and a subgroup of $\sum \ZZ e_b$. Hence $g\in Q$ 
acts on $H_1(P,\ZZ)$ via diagonal block matrices $diag(A(g),B(g))$, using the
intersection pairing one has that $B(g)={}^tA(g)^{-1}$ and $A(g)\in GL(4,\ZZ)$.
Therefore $H_1(P,\ZZ)\cong M\oplus M$.

In case $g(C)=g>2$ we obtain $C$ as a genus two surface, containing the loops $\alpha_{g-1}$ and $\alpha_g$,  with $g-2$ tori $T_1$, $\ldots$, $T_{g-2}$ attached. The inverse image of $T_i$ in $\tC$ is the disjoint union of 
of $8$ copies of $T_i$ which are permuted under the action of $Q$.
From this one finds the description given in the proposition.
\qed

\section{The Prym map is surjective in genus three.}\label{prym surj}

\subsection{} A quaternionic cover of a genus three curve produces an 
$8$-dimensional quaternionic Prym. These Pryms are parametrized by 
a $6$-dimensional symmetrical subdomain of the Siegel space, 
see section \ref{moduli}.
We will show in this section that the locus of such quaternionic Pryms 
has dimension $6$, hence the general principally polarized quaternionic 
$8$-fold in this family is a quaternionic Prym. 
In section \ref{diff} we show that it suffices to prove that a certain 
multiplication map $\mu_2$ is surjective and that is established in 
Theorem \ref{1.1}.

\subsection{}\label{diff}
 Let $\cA_g$ be the moduli space of principally polarized abelian 
varieties of dimension $g$. We will show that the map
$$
q:\cM_{3,Q}\longrightarrow \cA_8,\qquad [\tilde{C}\rightarrow C]
\longmapsto P
$$
has a 6-dimensional image. For this it suffices to find a cover
$\tilde{C}/C$ at which the differential ${\rm d}q$ of $q$ has maximal rank.
The map $q$ factors as follows:
$$
\cM_{3,Q}\stackrel{r}{\longrightarrow} \cM_{9}^\sharp\stackrel{p}{\longrightarrow} \cA_8,\qquad
 [\tilde{C}\rightarrow C]\longmapsto (\hC,\eta)\longmapsto P
$$
where $\cM_{9}^\sharp$ is the moduli space of genus 9 curves with a point of 
order two in the Picard group and $\eta\in Pic^0(\hC)[2]$
is the line bundle defining the 2:1 cover $\tilde{C}/\hC$.

The map $r$ is quasi-finite onto its image:
the curve $\hC$ comes with a fixed point free $(\ZZ/2\ZZ)^2$-action with 
quotient $C$,
as the automorphism group of a curve of genus at least 2 is finite, 
a given $\hC$ is in the image of at most a finite number of points in 
$\cM_{3,Q}$.

It remains to show that the differential of the Prym map $p$ has maximal rank
at a point $c=(\hC,\eta)$ in the image of $r$.

After adding level structures if necessary, the codifferential of $p$ 
$$
({\rm d}p_c)^*:T_{p(c)}^*\cA_8\cong S^2H^0(\omega_{\hC}\otimes \eta)\longrightarrow T_c^*\cM_{9}^\sharp\cong H^0(\hC,\omega_{\hC}^{\otimes 2})
$$
is the multiplication map $\mu_2:S^2H^0(\omega_{\hC}\otimes \eta)\rightarrow
H^0(\omega_\hC^{\otimes 2})$. The sheaf $\omega_{\hC}\otimes \eta$ is called
the Prym-canonical line bundle.
So we need to show that $\mu_2$ is surjective for
some point in the image of $r$.

\subsection{Lemma.}\label{spco}
There exist quaternionic covers $\tilde{C}/C$ in  $\cM_{3,Q}$ such that:
\begin{enumerate}
\item 
The genus $9$ curve $\hC=\tilde{C}/\{\pm 1\}$ is the complete intersection
of a smooth quadric $S$ and a quartic surface in $\PP^3$.

\item Let $L$ be the line bundle on $\hC$ which defines this map 
$\hC\rightarrow \PP^3$. Then $L$ is a theta characteristic 
on $\hC$ (so $L^{\otimes 2}\cong\omega$) and $h^0(L)=4$.

\item
Let $\theta:=L\otimes \eta$ where $\eta\in Pic^0(\hC)[2]$
is the line bundle defining the 2:1 cover $\tC/\hC$.
Then $\theta$ is a theta characteristic 
on $\hC$ and  $h^0(\theta)=0$.

\item
The Prym canonical line bundle $\omega_{\hC}\otimes \eta\cong
L\otimes \theta$ is very ample on $\hC$.
\end{enumerate}

\ts
Consider the genus two curve $X$ defined by $y^2=x^5-x$. We have $Q\subset Aut(X)$, for example:
$$
i,\,j:X\longrightarrow X,\qquad i(x,y):=(-x,iy),\quad j(x,y):=(x^{-1},iyx^{-3}).
$$
We will show that there are curves $\tC\subset J(X)$ which are invariant 
under $Q$ such that the quotient map 
$\tC\rightarrow C:=\tC/Q$ is unramified and such that $C$ has genus three. 
Thus the cover $\tC/C$ is in the family under consideration and we  verify the
lemma for these curves.
For such a curve $\tC$ the curve $\hC=\tC/\{\pm 1\}$ lies in the Kummer surface
$K(X)=J(X)/\pm 1$. If $\Theta\subset J(X)$ is (any) symmetrical theta divisor, the map $K(X)\subset\PP^3$  
given by $H^0(J(X),\cO(2\Theta))$ identifies $K(X)$ with a (16-nodal)
quartic surface.

The Kummer surface $K(X)\subset\PP^3$ is also the image of the secant variety 
$Sec(X)$ to the tricanonical curve $X\subset \PP^4$ under 
the map given by the quadrics in the ideal of $X$, this follows from \cite{BV},
Lemma 4.22(1) (this part of the lemma works for hyperelliptic curves),
see also \cite{vGI}, section 5 (with $D=K$, the canonical system) for another
proof.
Using the isomorphism 
$$
\langle 1,x,x^2,x^3,y\rangle\cong H^0(\omega_X^{\otimes 3}),
\qquad f\longmapsto fy^{-3}({\rm d}x)^{\otimes 3},
$$
one finds that $Q$ acts on $\PP^4$ via:
$$
i(x_0:\ldots:x_4)=(x_0:-x_1:x_2:-x_3:ix_4),\qquad
j(x_0:\ldots:x_4)=(x_3:x_2:x_1:x_0:ix_4).
$$
The ideal of the tricanonical curve $X\subset\PP^4$
is generated by the 4 quadrics
$$
Q_0=x_4^2+x_0x_1-x_2x_3,\quad Q_1=x_0x_2-x_1^2,\quad Q_2^2=x_0x_3-x_1x_2,
\quad Q_3=x_1x_3-x_2^2
$$
thus the rational map $\PP^4\rightarrow \PP^3$ given by $x\mapsto(Q_0(x):\ldots:Q_4(x))$ maps $Sec(X)$ onto $K(X)$.

The degree four polynomials $Q_0^2$, $Q_2^2$, $Q_1Q_3$ and $Q_1^2+Q_3^2$
are invariant under $Q$ and their common zeroes are just the points of 
$X\subset \PP^4$, in particular, they are not all zero on a secant line. 
Therefore a general linear combination of these quartic polynomials 
defines a smooth, irreducible $Q$-invariant divisor $\tC$ on $J(X)$
which does not contain any fixed point of any non trivial element of $Q$.
Hence $\tC\rightarrow C:=\tC/Q$ is a quaternionic cover.

Since $\tC \in |4\Theta|$,
the adjunction formula shows that $\omega_{\tC}$
is the restriction of $\cO_{J(X)}(4\Theta)$ to $\tC$ and thus
$\tC$ has genus $17$.
Moreover, $\tC$ has an even theta characteristic $\tilde{L}$ which is the 
restriction of 
$\cO_{J(X)}(2\Theta)$ to $\tC$. The exact sequence 
$$
0\longrightarrow \cO_{J(X)}(-2\Theta)\longrightarrow
\cO_{J(X)}(2\Theta)\longrightarrow \tilde{L}\longrightarrow 0
$$
shows that $h^0(\tilde{L})=h^0(\cO(2\Theta))=4$.
Since $2\Theta$ defines the map $J(X)\rightarrow K(X)\subset\PP^3$
which induces the cover $\pi:\tC\rightarrow \hC$, we have $\tilde{L}=\pi^*L$
for a theta characteristic $L$ on $\hC$ with $h^0(L)=4$.
As
$$
H^0(\tC,\tilde{L})\cong H^0(\hC,L)\oplus H^0(\hC,L\otimes\eta),
$$
we get $h^0(\theta)=0$ since $\theta=L\otimes \eta$.

To see that $L\otimes\theta$ is very ample, consider the minimal desingularization $\sigma: K \to K(X)$ of $K(X)$. Let $E$ be the exceptional divisor and let $D$ be the pull-back of a  plane. It turns out that
$$
\cO_K(4D-\mbox{$\frac 12$} E)
$$
is very ample and that its restriction to the pull-back of $C$ is $H$ cf.\ \cite{Ba}.
\qed

\subsection{Remark.} One can show that a genus $9$ curve $\hC$ 
which is a $(\ZZ/2\ZZ)^2$-cover of a genus $3$ curve $C$
has a theta characteristic $L$ with $h^0(L)=4$. 
The bundle $L$ is the pull-back of an odd theta characteristic on $C$.
The map $\hC\rightarrow |L|^*$
embeds the general such $\hC$ as a complete intersection of a smooth quadric and a quartic surface.

\subsection{}
From now on we write $C$ for the genus $9$ curve $\hC$ as in Lemma \ref{spco}.
We denote the Prym-canonical line bundle on $C$ by $H$:
$$
H=\omega_C\otimes \eta \cong L \otimes \theta.
$$
Since $H\not\cong \omega_C$, we have $h^1(H)=0$, hence $H$ is non special.
Moreover, $H^0(H\otimes L^*)=H^0(\theta)=0$. 
It follows from the previous lemma that   
$$
C \subset S = \PP^1 \times \PP^1,\qquad L=\cO_C(1,1),
$$
where $S$ is a smooth quadric and $C$ has type $(4,4)$.
The surjectivity of $\mu_2$ will be a consequence, 
for $n = 4$, of the following more general
statement about curves of type $(n,n)$ in $\PP^1 \times \PP^1$.

\subsection{Theorem.} \label {1.1} 
Let $C$ be any smooth element of $\mid \cO_S(n,n) \mid$ 
and let $L = \cO_C(1,1)$. Assume
$$
H \in Pic^{n^2}(C)
$$
is a non special and very ample line bundle such that $h^0(H \otimes L^*) = 0$. 

Then for all $d\geq 1$ the multiplication maps 
$$
\mu_d: Sym^d H^0(H) \to H^0(H^{\otimes d}) 
$$
are surjective, i.e.\ $H$ is normally generated.  

\ts
Since $C$ has genus $(n-1)^2$ and $H$ has degree $n^2$ with $h^1(H)=0$, Riemann-Roch shows that $h^0(H) = 2n$. As $H$ is very ample, we can assume
$$
C \subset \PP^{2n-1} 
$$  
with $\cO_C(1) \cong H$. Let $\cI_C$ be the ideal sheaf of $C$ in $\PP^{2n-1}$. It is easy to check that 
$$
h^0(\cI_C(2))\geq \dim S^2H^0(H)-\dim H^0(H^{\otimes 2}) = 2 \binom n2
$$ 
and that the equality holds iff $\mu_2$ is surjective.

Let $L_1 = \cO_C(1,0)$, $L_2 = \cO_C(0,1)$ then $L=L_1\otimes L_2$. 
We observe that, for any divisor
$$
d_i \in \mid L_i \mid \quad (i = 1,2),
$$
the linear span
$$
<d_i>\qquad(\subset\PP^{2n-1})
$$
is a space of dimension $n-1$. 
If not, $d_1 + d_2$ would be contained in a hyperplane $P$ and
we would have $P \cdot C = d_1+d_2+t$ with $t$ effective. 
This contradicts our assumptions,
because $\cO_C(t) \cong H \otimes L^* $. 
Using the same proof we can show that
$$
<d_1> \cap <d_2> 
$$
is either empty or a point of $C$. 
Assume $Supp d_1 \cap Supp d_2$ is empty, then the proof is exactly the same. 
Assume $Supp d_1 \cap Supp d_2$ is non empty, then it consists of one point
$p$ and moreover $Supp (d_1-p) \cap Supp (d_2-p)$ is empty. 
If $<d_1> \cap <d_2> = p$ we are done. 
Otherwise there exists a pencil of hyperplanes containing $d_1+d_2-p$. 
Hence there exists a hyperplane $P$ as above: contradiction. 

Let
$$
\Sigma_i
$$
be the union of the spaces $<d_i>$, where $d_i$ varies in $\mid L_i \mid$, then 
$$
C = \Sigma_1 \cap \Sigma_2.
$$
The variety $\Sigma_i$ is well known: it is an example of rational normal scroll of $\PP^{2n-1}$. 
A rational normal scroll $V$ of $\PP^r$ can be defined as the image, under the tautological map, of
$$
Proj \oplus_{s=1}^k \cO_{\PP^1}(a_s), 
$$
where $a_1 + \dots + a_k = r+1-k$ and $a_s \geq 0$. In particular $V$ has dimension $k$
and minimal degree $r+1-k$ with respect to the dimension. It is known that $V$ is
arithmetically Cohen-Macauley and that its homogeneous ideal
$$
I_V \subset B = \CC[X_0, \dots, X_r]
$$
is generated by quadrics. Moreover $V$ is a determinantal variety and a basis for
the vector space $I_V(2)$ is given by the order two minors, denoted by
$$
q_1, \dots ,q_N,
$$
of a $2 \times k$ matrix of linear forms (cfr. \cite{ACGH} p.95-100). In the case of $\Sigma_i$ 
we have $k = n$  so that $deg \Sigma_i = dim \Sigma_i = n$. Note that:
\begin{enumerate}
\item[(i)] $C$ is smooth of degree $n^2$, 
\item[(ii)] $codim C = codim \Sigma_1 + codim \Sigma_2$, 
\item[(iii)] $C = \Sigma_1 \cap \Sigma_2$.
\end{enumerate}
 This implies that $\Sigma_1$ and $\Sigma_2$ intersect
transversally at each point and that $C$ is scheme theoretically 
the intersection of $\Sigma_1$ and $\Sigma_2$. 
Moreover the depth of $C$, as a subscheme of $\Sigma_1$, 
equals its  codimension in $\Sigma_2$, which is $n-1$. 
(We recall that, by definiton, the depth of a closed subscheme $Y$ in $X$ 
is the largest integer 
$d$ such that, for all $x \in Y$, the ideal of $Y$ in the local ring  
${\cal O}_{X,x}$ contains a regular sequence of length $d$: 
see \cite{Fulton} 14.3.1, p.251). 

The ideal $I_{\Sigma_2}$ admits a minimal free resolution 
$(F_{*})Ê\to I_{\Sigma_2}$
of the same length $n-1$. As is well known such a resolution looks like
$$
0 \to F_{n-1} \to \dots \to F_1 \to I_{\Sigma_2} \to 0,
$$ 
where $F_1$ is the direct sum of $\binom n2$ copies of $B(-2)$ and 
the map $F_1 \to I_{\Sigma_2}$
is defined by the previously considered minors $q_1, \dots, q_N$. 
Let
$$
J \subset A = B/I_{\Sigma_1}
$$
be the homogeneous ideal of $C$ in $\Sigma_1$. 
We claim that the tensor product $(F_{*})\otimes_B A$ yields 
a resolution of $J$. 
This is true because the length of $(F_{*})$ and the depth of $C$ in 
$\Sigma_1$ are the same. Indeed, under this assumption, we are in 
the case considered in \cite{Fulton}, Example 14.3.1 p.252. 
Our claim follows from the property stated there or 
from \cite{KL}  Corollary 8. 

Let $I_C$ be the homogeneous ideal of $C$ in $\PP^{2n-1}$ and let $q \in I_C(2)$. Due to our
previous claim, $q$ is a linear combination of elements of 
$I_1(2)$ and $I_2(2)$. Hence we have $dim I_C(2) \leq 2 \binom n2$. 
On the other hand we had $dim I_C(2) \geq 2  \binom n2$ so
that the equality holds. Hence $\mu_2$ is surjective and this implies that
$\mu_d$ is surjective for all $d\geq 2$ (cfr. \cite{ACGH} D.6 p.140).  
\qed

\subsection{Remark.} Let $g$ be the genus of $C$, 
then $g = (n-1)^2$ and the line bundle $H \otimes L^*$
lies in $Pic^{g-1}(C)$. It can be shown that the assumptions 
of the theorem are satisfied
on a non empty open set of $Pic^{g-1}(C)$. \par \noindent
Conversely, every transversal intersection of two rational normal scroll of dimension $n$ 
in $\PP^{2n-1}$ admits a natural embedding in $\PP^1 \times \PP^1$ as a smooth curve of type $(n,n)$.

\section{Weil classes and quaternion classes}\label{weil}

\subsection{} We recall a result of Abdulali \cite{A} which relates 
the Hodge conjecture
 for abelian varieties of Weil type to the Hodge conjecture for 
 quaternionic abelian varieties. 
From Schoen's construction of cycles on Prym varieties and 
our previous results, the Hodge conjecture for the general quaternionic Prym 
studied in this paper then follows, see Corollary \ref{cycler}.
 
\subsection{Abelian varieties of Weil type.}
For an imaginary quadratic field $K$ we denote by $x\mapsto\bar{x}$ the 
non-trivial field automorphism of $K$. The complex embeddings $K\subset\CC$ are denoted by $\sigma$ and $\bar{\sigma}$, note $\sigma(\bar{x})=\bar{\sigma}(x)$.

A polarized abelian variety $(A,E)$
of dimension $2n$ is of Weil type, with field $K$,
 if there is an inclusion $K\subset End(A)\otimes \QQ$ such that $H_1(A,\QQ)$ is a $K$-vector space and if the eigenvalues of any $x\in K$ on the holomorphic tangent space $T_0A$ to $A$ at the origin are $\sigma(x)$ and $\bar{\sigma}(x)$
each with multiplicity $n$. Moreover, one requires that $x^*E=x\bar{x} E$.

\subsection{Weil classes.}
We recall the definition of the space of Weil classes $W_K\subset H^{2n}(A,\QQ)$
for an abelian variety $A$ of Weil type of dimension $2n$ with field $K$. 

For an embedding $\rho$ of $K$ into $\CC$ we let $V_\rho\subset H^1(A,\CC)$ be the corresponding eigenspace, so $x\cdot v=\rho(x)v$ for $x\in K$ and $v\in V_\rho$:
$$
H^1(A,\QQ)\otimes_\QQ\CC\cong H^1(A,\CC)=
\oplus_{\rho\in \{\sigma,\bar{\sigma}\}} V_\rho,\qquad \dim_\CC V_\rho=2n.
$$
Since $A$ is of Weil type, the Hodge decomposition induces a decomposition
$$
V_\rho=V_\rho^{1,0}\oplus V_\rho^{0,1},\qquad 
\dim_\CC V_\rho^{1,0}=\dim_\CC V_\rho^{0,1}=n
\qquad(\rho\in \{\sigma,\,\bar{\sigma}\}).
$$ 

It is easy to see that the one dimensional subspaces $\wedge^{2n}V_\sigma$ and
$\wedge^{2n}V_{\bar{\sigma}}$ of $H^{2n}(A,\CC)$ are of Hodge type $(n,n)$
and that they are the eigenspaces in $H^{2n}(A,\CC)$ on which $x\in K$ acts
as $\sigma(x)^{2n}$ and $\bar{\sigma}(x)^{2n}$.
This implies that their sum is defined over $\QQ$. In fact,
for $x\in K$, the polynomial $f:=(T-\sigma(x)^{2n})(T-\bar{\sigma}(x)^{2n})$ 
has coefficients in $\QQ$, let $f=T^2+aT+b$. 
The kernel $W_x$ of the linear map
$$
W_x=\ker\left(x^2+ax+b: H^{2n}(A,\QQ)=\wedge^{2n}H^1(A,\QQ)\longrightarrow \wedge^{2n}H^1(A,\QQ)\right)
$$
(induced by $x^2+ax+b\in End(H^1(A,\QQ))$) is the subspace of $\wedge^{2n}H^1(A,\QQ)$
on which the eigenvalues of $x$ are $\sigma(x)^{2n}$ and $\bar{\sigma}(x)^{2n}$.
Therefore $W_K$, the intersection of all $W_x$ ($x\in K$), is two dimensional:
$$
W_K:=\cap_{x\in K} W_x,\qquad W_K\otimes_\QQ\CC=
\; \wedge^{2n}V_\sigma\,\oplus \wedge^{2n}V_ {\bar{\sigma}}.
$$
The subspace $W_K$ is also a sub-Hodge structure since each $x\in K$ preserves the Hodge structure on $H^1(A,\QQ)$ and, as we observed before:
$$
W_K\otimes\CC=\oplus_{\rho\in\{\sigma,\,\bar{\sigma}\}}\wedge^{2n}V_\rho=
\oplus_{\rho\in\{\sigma,\,\bar{\sigma}\}}
\wedge^{n}V_\rho^{1,0}\otimes \wedge^{n}V_\rho^{0,1}\subset H^{n,n}(A),
$$
hence $W_K\subset H^{n,n}(A)$.
Thus we have a two dimensional subspace $W_K\subset H^{2n}(A,\QQ)\cap H^{n,n}(A)$, which is the space of Weil classes. 
(To see that this is the space defined by Weil in \cite{We}, note that the subspace $\wedge^{2n}_K H^1(A,\QQ)$ is contained in $H^{2n}(A,\QQ)$ and
that $x\in K$ has eigenvalues $\rho(x)^{2n}$ on this subspace.)

\subsection{Hodge classes.}
For a general $2n$-dimensional polarized abelian variety $(A,E)$ of Weil type 
(i.e.\ $Hod(A)(\CC)\cong SL(2n,\CC)$) one has (\cite{vG}, Thm.\ 6.12):
$$
B^n(A)\cong \langle \wedge^n E\rangle \,\oplus\, W_K,\qquad
B^i(A)=\langle \wedge^iE\rangle\qquad(i\neq n).
$$

\subsection{Lemma.} Let $(A,E,F)$ be an abelian variety of quaternion type.
Let $K\subset F$ be a quadratic extension of $\QQ$. Then $A$ is of Weil type for
$K$.

\ts
A more general statement is proved in \cite{A}, $\S$4.
Since $F$ is definite, $K\subset F$ must be an imaginary quadratic field.
Moreover, there is a $j\in F$ such that $F=K\oplus Kj$ and $xj=j\bar{x}$
for $x\in K$. Hence the eigenspaces of $x$ are permuted by the action of $j$ on $T_0A$ and therefore they have the same dimension.
\qed

\subsection{Quaternion classes.}
Let $F$ be a definite quaternion algebra over $\QQ$ and let $K\subset F$ be an (imaginary) quadratic subfield. 
Let $A$ be a $2n$-dimensional abelian variety of quaternion type with algebra $F$. We define the subspace $W_F\subset H^{2n}(A,\QQ)$ as the subspace spanned by the translates $x\cdot W_K$ where $x$ runs over $F$ and $W_K$ is the space of Weil classes for the field $K$. 

\subsection{Proposition.}
The definition of $W_F$ does not depend on the choice of $K\subset F$ and
$$
W_F\subset H^{2n}(A,\QQ)\cap H^{n,n}(A),\qquad\dim_\QQ W_F=2n+1.
$$

\ts
Since $\HH_\QQ\otimes_\QQ\CC\cong M_2(\CC)$ we have $(F\otimes\CC)^*\cong GL_2(\CC)$ and this group acts on each $\wedge^kH^1(A,\CC)$. The space $H^1(A,\CC)\cong V^{2n}$ where $V$ is the standard 2-dimensional representation
of $GL_2(\CC)$. If $v_i$ denotes a highest weight vector in the $i$-th copy of $V$ then $v_1\wedge\ldots\wedge v_{2n}$ is a highest weight vector in
$\wedge^{2n} V^{2n}$ and the multiplicity of the corresponding weight is one.
Hence $\wedge^{2n}H^1(A,\CC)$ has a unique irreducible summand 
of dimension $2n+1$ denoted by $V_{2n+1}$.

Let $K\subset F$ be a quadratic subfield, then $(K\otimes \CC)^*\cong(\CC^*)^2$
and this is a Cartan subgroup of $GL_2(\CC)$. Since $x\in K$ has eigenvalues
$\sigma(x),\,\bar{\sigma}(x)$ on $V$, it has eigenvalues $\sigma(x)^{2n}$,
$\sigma(x)^{2n-1}\bar{\sigma}(x)$,\ldots, $\bar{\sigma}(x)^{2n}$ on $V_{2n+1}$.
In particular, $W_K\otimes \CC\subset V_{2n+1}$.
The uniqueness and irreducibility of $V_{2n+1}$ imply that $W_F\otimes\CC= V_{2n+1}$, hence $W_F$ does not depend on $K$ and has dimension $2n+1$. That
$W_F\subset H^{n,n}(A)$ is obvious since $W_K\subset H^{n,n}(A)$ and
elements from $F$ preserve the Hodge structure on $H^{2n}(A,\QQ)$.
\qed

\subsection{Quaternion Hodge classes.}
For a general polarized  $2n$-dimensional abelian variety $(A,E)$ 
of quaternion type 
(i.e.\ $Hod(A)(\CC)\cong SO(2n,\CC))$ one has (\cite{A}, Thm 4.1):
$$
B^n(A)\cong \langle \wedge^n E\rangle \,\oplus\, W_F,\qquad
B^i(A)=\langle \wedge^iE\rangle\qquad(i\neq n).
$$

\subsection{Corollary.}\label{wq}(Abdulali \cite{A}.)
Let $(A,E,F)$ be an abelian variety of quaternion type.
Assume there is a quadratic subfield $K\subset F$ such that
the space of Weil classes $W_K$ is spanned by cycle classes.

Then the space of quaternionic classes $W_F$ is spanned by cycle classes.

\ts This follows from the fact that $W_F$ is spanned by the translates
of $W_K$ by elements of $F\subset End(A)$.
\qed

\subsection{Corollary.}\label{cycler}
Let $\tC\rightarrow C$ be a quaternionic cover of a genus $g$ curve $C$ and let
$(P,E,F=\HH_\QQ)$ be the associated quaternionic Prym as in section \ref{qp}.

Then the subspace $W_F\subset H^{4(g-1)}(P,\QQ)$ is spanned by cycle classes.

\ts 
Let $K=\QQ(i)\subset \HH_\QQ$, and consider the 4:1 unramified cover $\pi_i:\tilde{C}\rightarrow \tilde{C}/\langle i\rangle$.
By Schoen's result \cite{S1}, applied to the cover $\pi_i$, the space of Weil classes $W_K\subset H^{4(g-1)}(P,\QQ)$ is spanned by cycle classes. 
Now apply Corollary \ref{wq}.
\qed

\section{Specialisation of quaternion classes}

\subsection{} 
We study the space of quaternionic classes $W_F\subset H^{2n}(A)$ 
in the case that the quaternionic variety $A$ is 
a self-product $B^2$ of an abelian $n$-fold of Weil type $B$. 
The result is disappointing: the classes in 
$W_F$ specialize to complete intersections of divisors on $B^2$. 
In particular, even if the classes in $W_F$ are known to be cycle classes, 
one does not obtain new cycle classes on $B^2$.

\subsection{}
Let $K\subset F$ be an imaginary quadratic subfield of a definite quaternion algebra $F$ over $\QQ$.
Let $A$ an abelian variety of quaternion type with algebra $F$
which is isogenous to $B^2$, where $B$ is of Weil type with $End(B)_\QQ\cong K$. Then we have 
$$
End(A)_\QQ=M_2(K)\;\cong F\otimes_\QQ K.
$$

The isomorphism between the algebras can be obtained as follows. 
Let $K\cong \QQ(\sqrt{r})$, and write $F=K\oplus Kj$ as in \ref{defs}, 
with $j^2=s$ and $jx=\bar{x}j$. The map
$$
F\,\hookrightarrow \,M_2(K),\qquad
x+yj\longmapsto \left(
\begin{array}{rcl}
x&y\\ s\bar{y}&\bar{x}
\end{array}
\right)
$$
is an injective $\QQ$-algebra homomorphism and identifies $F$ with a subalgebra
of $M_2(K)$. 
The $K$-linear extension of this map defines an isomorphism 
$F\otimes_\QQ K\cong M_2(K)$. 

The idempotents
$$
P_1:= \left(
\begin{array}{rcl} 1&0\\0&0\end{array}
\right)\qquad{\rm and}\quad P_2:=I-P_1
$$
define, up to isogeny, subvarieties $B_i=Im(P_i)$ of $A$ and $A$ is isogenous to $B_1\times B_2$. The uniqueness of the decomposition of $A$ into simple factors implies that the abelian varieties $B_1$ and $B_2$ are isogenous to the simple variety $B$.
An explicit isogeny 
between $B_1$ and $B_2$ is given by $j\in F$. 
Note that for $x\in K$ the action of
$x\otimes 1$ and $1\otimes x$, coincide on $B_1$ but that on $B_2$ the action of $x\otimes 1$ and $1\otimes \bar{x}$ coincide.

\subsection{Proposition.} \label{so8res}
Let $(A,F)$ be an abelian variety of quaternion type
of dimension $2n$, let $K$ be a quadratic subfield of $F$ and assume that 
$A$ is isogenous to $B^2$ with $B$ of Weil type for the field $K$.
 
Then the subspace $W_F\subset H^{2n}(A,\QQ)$ is spanned by intersections of divisor classes.

\ts 
To determine the divisor
classes in $H^2(A,K)$ we use
that $(H^2(A,\QQ)\cap H^{1,1}(A))\otimes\CC$ is the subspace of 
the invariants of the Hodge group $Hod(A)(\CC)$ in $\wedge^2H^1(A,\CC)$:
$$
(H^2(A,\QQ)\cap H^{1,1}(A))\otimes\CC=\wedge^2H^1(A,\CC)^{Hod(A)(\CC)}.
$$
Since $A\sim B^2$, the Hodge group $Hod(A)(\CC)$ is isomorphic to
$Hod(B)(\CC)$, acting diagonally on $H^1(A,\CC)=H^1(B,\CC)\oplus H^1(B,\CC)$.
As $B$ is of Weil type, $Hod(B)(\CC)\subset SL(n,\CC)$ 
 and $H^1(B,\CC)$ is isomorphic to the direct sum of the standard representation $W$ of $SL(n,\CC)$ and its dual $W^*$ (cf.\ \cite{vG}). 
The $SL(n,\CC)$-invariants in $H^2(A,\CC)$ are thus also $Hod(A)(\CC)$-invariants.
We have:
$$
\wedge^2H^1(A,\CC)\cong (\wedge^2 W)^2\oplus (\wedge^2 W^*)^2\oplus
(W\otimes W)\oplus (W^*\otimes W^*)\oplus (W\otimes W^*)^4.
$$
There are no $SL(n,\CC)$-invariants in $\wedge^2 W$, $\wedge^2 W^*$, $W\otimes W$ and $W^*\otimes W^*$, but the invariants are one dimensional in $W\otimes W^*$
(identifying this space with $End(W)$, these are the scalar multiples of the identity). Thus $\dim H^2(A,\QQ)\cap H^{1,1}(A)\geq 4$.

Next we show that there is a one dimensional subspace in 
$ H^2(A,\CC)\cap H^{1,1}(A)$ on which the subfield $K\subset F$ acts 
via $x\cdot v=x^2v$.
Under the isogeny $A\sim B_1\times B_2$ given by the projectors $P_i$
we have
$$
H^1(B_1,\QQ)\cong V_1:=\{v\in H^1(A,\QQ):\;(x\otimes 1)v=(1\otimes x)v,\quad
\forall x\in K\}, 
$$
$$
H^1(B_2,\QQ)\cong V_2:=\{v\in H^1(A,\QQ):\;(x\otimes 1)v=(1\otimes\bar{x})v,\quad
\forall x\in K\}.
$$
The subspaces $H^1(B_i,\CC)=V_i\otimes \CC$ of $H^1(A,\CC)$ have an eigenspace decomposition:
$$
V_i\otimes_\QQ \CC=V_i^+\oplus V_i^-,\qquad (x\otimes 1)v_i^+=xv_i^+,\quad
(x\otimes 1)v_i^-=\bar{x}v_i^-,
$$
for all $v_i^\pm\in V_i^\pm$ and $i=1,\,2$. These subspaces are subrepresentations of $Hod(A)(\CC)$ (since the Hodge group commutes with the
endomorphisms) and we may identify $V_1^+\cong W$, $V_1^-\cong W^*$ because $B_1$ is of Weil type.
The problem is to find the sign for which $V_2^\pm\cong W$. Using the action
of $j\in F$ this is easy though. Since $j\in End(A)$ it commutes with the action of the Hodge group, also $jV_1^+\cong W$. As $xj=j\bar{x}$ in $F$, we have $jV_1^+= V_2^-$. Therefore $V_2^-\cong W$ and $V_2^+\cong W^*$.

This implies that there is a one-dimensional subspace of $V_1^+\otimes V_2^+$ which is contained in $H^2(A,\CC)\cap H^{1,1}(A)$. The action of 
$K\subset F$ on the vector space $V_1^+\otimes V_2^+$ is via $x\cdot v=x^2 v$,
as desired. 

Finally we let $c$ be a basis of this one dimensional vector space.
The $n$-fold self-product of $c$ is $c^n\in \wedge^{2n}H^1(A,\CC)\cap H^{n,n}(A)$. The field $K$ acts on $c^n$ as $x\cdot c^n=x^{2n}c^n$, thus
$c^n\in W_K\otimes\CC$. 
Since $c\in V_1^+\otimes V_2^+\cong W\otimes W^*$ corresponds to the identity map, we have that $c^n\in \wedge^n V_1^+\otimes \wedge^n V_2^+$ $\subset \wedge^{2n}(V_1^+\oplus V_2^+)$ is non-trivial 
and lies in $W_K$.  
Hence $W_K$ is spanned by divisor classes.
Since $W_F$ is spanned by $F$-translates of $W_K$, we conclude that $W_F$ is spanned by divisor classes.
\qed

\section{Abelian 8-folds of so(7)-type}

 \subsection{} 
The general abelian 8-fold $A$ of quaternion type has complex Hodge group 
$Hod(A)(\CC)\cong SO(8,\CC)$ and the representation of this group on 
$H^1(A,\CC)$ can be identified, using triality, 
with 2 copies of a half spin representation $\Gamma_+$ of dimension $8$. 
In the standard representation of $SO(8,\CC)$ on $\CC^8$ there is an 
invariant quadratic form and the stabilizer of a(ny) non-isotropic vector 
in $\CC^8$
is a subgroup isomorphic to $SO(7,\CC)$. 
The half spin representation $\Gamma_+$ of $SO(8)$ restricts to the spin 
representation $\Gamma$ of $SO(7)$,
an eight dimensional irreducible representation. 
Now it is possible to find 5-dimensional families of abelian 8-folds $A$
whose complex Hodge group is $SO(7)$ and its representation on $H^1(A,\CC)$
is two copies of the spin representation of $SO(7)$.  

We will say that an $8$-dimensional abelian variety $A$ is of $so(7)$ type 
if the Lie algebra of $Hod(A)(\CC)$ is contained in $so(7)$ and if $H^1(A,\CC)$ is isomorphic to two copies of the spin representation of $so(7)$.

\subsection{Proposition.}
In a $6$-dimensional family of 8-folds of quaternion type there
exist $5$-dimensional families of abelian varieties of $so(7)$ type.
Moreover, in such a $5$-dimensional family there are $4$-dimensional families 
of abelian varieties of Weil type.

\ts
We use the theory of Kuga-Satake varieties but alternatively one could directly use Mumford-Tate groups. The moduli space of polarized Hodge structures of weight $2$ with Hodge numbers $(1,n,1)$ has dimension $n$ and the complex Hodge group of a general such Hodge structure is $SO(n+2)$ (\cite{vG2}, 4.6, 4.7). 
Kuga and Satake associated to such a Hodge structure $V$ an abelian variety 
$J_{KS}(V)$ of dimension $2^{n-2}$ whose complex Hodge group is the universal 
cover $Spin(n+2)$ of $SO(n+2)$. 
This abelian variety is not simple in general, one has for $V$ general and
$n=4,\,5,\, 6$ respectively:
$$
J_{KS}(V)\sim B_4^4,\quad  A^4, \quad (A_+\times A_-)^4
$$
where $B_4$ is an abelian $4$-fold of Weil type (\cite{vG2}, Theorem 9.2),
$A$ is an abelian $8$-fold of $so(7)$ type and $A_+$, $A_-$ are non-isogenous abelian varieties of quaternion type with complex Hodge group isomorphic to $SO(8)$. 
Specializing a general Hodge structure of type $(1,6,1)$ to a direct sum of a
general Hodge structure of type $(1,7,1)$ and a Tate Hodge structure $\QQ(-1)$,
the abelian varieties $A$ and $A'$ become isogenous and of $so(7)$-type, similarly, specializing a general Hodge structure of type $(1,5,1)$ to one of type $(1,4,1)$ and $\QQ(-1)$, the Kuga-Satake variety specializes to $B_4^8$
where $B$ is of Weil type (note that $Hod(B)(\CC)=SL(4,\CC)=Spin(6,\CC)$). 
\qed

\subsection{}
The Hodge cycles $B^p(A)=H^{2p}(A,\QQ)\cap H^{p,p}(A)$ 
in codimension $p$ on $A$ are just the invariants of the Hodge group $Hod(A)$, 
hence:
$$
B^p(A)\otimes_\QQ\CC=H^{2p}(A,\CC)^{Hod(A)(\CC)}=
\left(\wedge^{2p} (\Gamma\oplus\Gamma) \right)^{so(7)},
$$
where $so(7)$ is the Lie algebra of $SO(7,\CC)$. 
The following proposition collects the basic facts on $so(7)$-representations
that we will need.

\subsection{Proposition.}
Let $\Gamma$ be the $8$-dimensional spin representation and
$V$ is the standard $7$-dimensional representation of $so(7)$. Then:
$$
\wedge^2\Gamma=\wedge^2 V\,\oplus V,\qquad S^2\Gamma=\wedge^3 V\,\oplus \,\CC
$$
is the decomposition of $\wedge^2\Gamma$ in irreducible 
$so(7)$-representations,
here $\CC$ is a trivial one dimensional representation.
Similarly:
$$
\wedge^4\Gamma\cong 
\wedge^3 V\;\oplus\;(S^2V)_0\;\oplus\; V \oplus\; \CC,
$$
where $S^2V=(S^2V)_0\oplus \langle Q\rangle$.  

\ts The decomposition of $\wedge^2\Gamma$ can be found as in
 \cite{FH}, sections 19.3 and 19.4.
The decomposition of $\wedge^4\Gamma$ follows from 
an easy explicit calculation which shows that both sides 
have the same weights and multiplicities.
\qed

\subsection{Corollary.}
The dimensions of the spaces of Hodge classes on $A$ are:
$$
\dim B^1(A)=1,\qquad \dim B^2(A)=6.
$$
Moreover, $End(A)$ is four dimensional algebra over $\QQ$.

\ts
Note that $H^2(A,\CC)$ is the $so(7)$-representation:
$$
\wedge^2 (\Gamma\oplus\Gamma)=\left(\wedge^2\Gamma\right)^2\oplus
(\Gamma\otimes\Gamma)=\left(\wedge^2\Gamma\right)^3\oplus Sym^2(\Gamma).
$$
There are no $so(7)$-invariants in $\wedge^2\Gamma$ and there is one invariant 
in $Sym^2\Gamma$ given by the quadratic form $Q_s\in Sym^2\Gamma$.
Hence the Picard number of $A$ is  $\dim B^1(A)=1$.  
Furthermore, $\dim End(A)=4$ (compute invariants in $H^1(A)\otimes H^1(A)$).

To determine $\dim B^2(A)$, we decompose $H^4(A,\CC)$:
$$
\wedge^4 (\Gamma\oplus\Gamma)=\left(\wedge^4\Gamma\right)^2\;\oplus\;
\left((\wedge^3\Gamma)\otimes\Gamma\right)^2 \;\oplus\; 
(\wedge^2\Gamma)\otimes(\wedge^2\Gamma).
$$
Since $\dim \left(\wedge^4\Gamma\right)^{so(7)}=1$ and
$\wedge^4\Gamma$ is a summand of $(\wedge^3\Gamma)\otimes\Gamma$
and of $(\wedge^2\Gamma)\otimes(\wedge^2\Gamma)$, we find $\dim B^2(A)\geq 5$.
One can verify that $\dim (\wedge^2\Gamma\otimes\wedge^2\Gamma)^{so(7)}=2$
and $\dim \left((\wedge^3\Gamma)\otimes\Gamma\right)^{so(7)}=1$,  
hence $\dim B^2(A)=6$. The program \cite{LiE} may be helpfull here.
\qed

\subsection{Remark.} One can show that $\dim B^3(A)=6$ and $\dim B^4(A)=16$
for an abelian $8$-fold of $so(7)$ type.

\subsection{}
Since $\dim B^1(A)=1$, we have $B^1(A)=\langle E\rangle$ where $E$ is a 
polarization of $A$.
The following lemma shows that if there exists an algebraic cycle on $A$ 
whose class $c$ is not in $\langle \wedge^2E\rangle \subset B^2$, 
then $B^2(A)$ is spanned by algebraic cycles. 
 
 \subsection{Lemma.}
Under the action of the multiplicative group $End(A)^*$ of the quaternion algebra $End(A)$, the space $B^2(A)$ has two irreducible components, one of dimension $1$ (spanned by $\wedge^2E$ where $B^1(A)=\langle E\rangle$) 
and the other is irreducible of dimension $5$ and will be denoted by 
$B^2(A)_0$:
$$
B^2(A)=\langle \wedge^2E\rangle\;\oplus\; B^2(A)_0.
$$

\ts
Note that $(End(A)\otimes\CC)^*\cong GL(2,\CC)$, compatible with the splitting
$H^1(A,\CC)=\Gamma\oplus\Gamma$. Since the diagonal matrix 
$h(t)=diag(t,t^{-1})$
acts as $t^4$ on $\wedge^4$ of the first factor, 
$GL(2,\CC)$ has an irreducible $5$-dimensional representation 
in $B^2(A)\otimes\CC$.
\qed

\subsection{Theorem.} \label{thm7}
Let $A$ be an abelian $8$-fold of $so(7)$-type and assume that $A$ is 
isogenous to $B^2$ where $B$ is of Weil type for an 
imaginary quadratic field $K\subset F=End(A)$.

If $B^2(A)$ is spanned by classes of algebraic cycles, 
then the Weil classes $W_K\subset B^2(B)$ are classes of algebraic cycles.

\ts
Let $E_B$ be a polarization on $B$. Then 
$$
B^2(B)=\langle \wedge^2E_B\rangle\oplus W_K
$$
and the eigenvalue of $x\in K^*=End(B)^*$ on $\wedge^2E_B$ is $(x\bar{x})^2$
and its eigenvalues on $W_K$ are $x^4$ and $\bar{x}^4$.
In particular, the two components of an algebraic cycle class in $B^2(B)$ are 
again algebraic cycle classes.
To prove the theorem it thus suffices to show that the image of the composition
$$
B^2(A)\hookrightarrow H^4(A,\QQ)\stackrel{\kappa}{\longrightarrow} 
H^4(B_1,\QQ)\otimes H^0(B_2,\QQ)\stackrel{i_1^*}{\longrightarrow }
H^4(B_1,\QQ)
$$
is not contained in $\langle \wedge^2E_1\rangle\subset B^2(B_1)$, 
here $\kappa$ is obtained from the K\"unneth decomposition 
$H^4(A)=\oplus H^i(B_1)\otimes H^{4-i}(B_2)$ 
(note that the K\"unneth components of an algebraic cycle on $B_1\times B_2$
are algebraic, use the endomorphisms 
$(b,b')\mapsto (nb,mb')$ with $n,\,m\in\ZZ$) and
$i_1:B_1\rightarrow B_1\times B_2$, $b\mapsto (b,0)$. 

Recall that $H^1(A,\CC)\cong \Gamma_1\oplus\Gamma_2$ with $\Gamma_i\cong \Gamma$, the spin representation fo $so(7)$,  and that
$\dim (\wedge^4\Gamma)^{so(7)}=1$. Let 
$$
\langle  F\rangle =\wedge^4\Gamma_1\hookrightarrow  B^2(A)\otimes \CC.
$$
It suffices to show that $F$ has a non-trivial component in $W_K\otimes\CC\subset H^4(B_1,\CC)$. As $Hod(B_1)(\CC)=SL(4,\CC)$ representations, we have $H^1(B_1,\CC)=W_1\oplus W_1^*$ and 
$W_K\otimes\CC=(\wedge^4W_1)\oplus(\wedge^4W_1^*)$.

We now indicate a somewhat subtle point: $\Gamma_1$ does not restrict to 
$W_1\oplus W_1^*$, but nevertheless, $W_1$ is an irreducible component of 
the restriction of $\Gamma_1$.

For a general $A$ of $so(7)$-type,
the $\QQ$-vector space $H^1(A,\QQ)$ is an irreducible $Hod(A)$-representation.
The actions of $Hod(A)$ and $F$ commute and
extend $\CC$-linearly to the vector space $H^1(A,\CC)=H^1(A,\QQ)\otimes\CC$. The eigenvalues of the action of $x\in K$ on $H^1(A,\QQ)$  
are $\sigma(x)$ and $\bar{\sigma}(x)$, hence we have the eigenspace 
decomposition:
$$
H^1(A,\CC)\cong V_\sigma\oplus V_{\bar{\sigma}} 
\qquad{\rm with}\quad
V_\rho:=\{v\in H^1(A,\CC):\;x\cdot v=\rho(x)v\;\}.
$$
Since the action of $K$ commutes with the action of $Hod(A)$, each of these eigenspaces is a $Hod(A)(\CC)$-representation and we will identify
$\Gamma_1=V_\sigma$.

In case $A$ is isogenous to $B^2$, we get $End(A)\cong M_2(K)\cong F\otimes K$ and the center
$K$ of $M_2(K)$ also acts on $A$. This leads to a further decomposition:
$$
H^1(A,\CC)\cong W_{\sigma,\sigma}\oplus W_{\sigma,\bar{\sigma}}
\oplus W_{\bar{\sigma},\sigma}\oplus W_{\bar{\sigma},\bar{\sigma}}
$$
where, with $x\otimes y\in K\otimes K\subset F\otimes K$:
$$
W_{\rho,\tau}:=\{v\in H^1(A,\CC):\; (x\otimes y)\cdot v=\rho(x)\tau(y)v\quad
x,\,y\in K\,\}.
$$
Note that $x\otimes 1$ and $1\otimes x$ act in the same way on both $W_{\sigma,\sigma}$ and $W_{\bar{\sigma},\bar{\sigma}}$, hence, 
with the notation of the proof of \ref{so8res}, we have:
$$
H^1(B_1,\QQ)=W_{\sigma,\sigma}\oplus W_{\bar{\sigma},\bar{\sigma}},\qquad
V_1^+=W_{\sigma,\sigma},\quad V_1^-=W_{\bar{\sigma},\bar{\sigma}}.
$$
Note that $W_{\bar{\sigma},\bar{\sigma}}$ is  {\em not} contained in $\Gamma_1=V_\sigma$.
Thus the $so(7)$-representation $\Gamma$ restricts to the 
$so(6)=sl(4)$-representation $W\oplus W^*$ and, although
this representation is isomorphic to $H^1(B,\CC)$,
we actually have:
$$
H^1(B_1,\CC)=W_1\oplus W_2^*,\qquad H^1(B_2,\CC)=W_1^*\oplus W_2
$$
where $\Gamma_i$ restricts to $W_i\oplus W_i^*$. In any case,
we have to verify that the $so(7)$-invariant $F\in \wedge^4\Gamma_1$ 
has a non-trivial component in $\wedge^4W_1$.

We use the explicit description of the spin representation given in \cite{FH}. 
Let $Q=x_7^2+\sum^3_{i=1} x_ix_{i+3}$ be the quadratic form on $\CC^7$ 
defining $so(7)$. 
Let $R=\langle e_1,e_2,e_3\rangle$,
then $R$ is an isotropic subspace of $Q$, and let
$$
S=\oplus_{i=0}^3 \wedge^i R\;\cong\;\CC^8.
$$ 
The spin representation of $so(7)$ is realized on $S$.
A basis of weight vectors of $S$ is given by the 
$$
e_I=e_{i_1}\wedge\ldots\wedge e_{i_k}\qquad{\rm if}\quad
I=\{i_1,\ldots,i_k\}\subset\{1,\,2,\,3\},\quad i_1<\ldots <i_k,
$$
in fact, $e_I$ has weight $(1/2)(\sum_{i\in I} L_i- \sum_{j\not\in I} L_j)$
(\cite{FH}, proof of Prop.\ 20.20).
Since an invariant has weight $0$, $F$ must lie in the weight $0$-subspace of 
$\wedge^4\Gamma$ which has dimension $8$.
 
Using the isomorphism $so(7)\cong \wedge^2\CC^7$, the spin representation 
$\Gamma$ is realized on $S$ by explicit formulas (cf.\ \cite{FH},
lemma 20.16), for example:
$$
(e_i\wedge e_7)\cdot e_I=(-1)^{\sharp I}e_i\wedge e_I\qquad (1\leq i\leq 3)
$$
where $\sharp I$ is the cardinality of $I$.
It is then very easy (but a bit tedious) to check that the only element
in the weight $0$-subspace which is mapped to zero by the $e_i\wedge e_7$
($1\leq i\leq 3$) are the scalar multiples of:
$$
2e_\emptyset\wedge e_{\{1,2\}}\wedge e_{\{1,3\}}\wedge e_{\{2,3\}}\!-
e_\emptyset\wedge e_{\{3\}}\wedge e_{\{1,2\}}\wedge e_{\{1,2,3\}}\!+
e_\emptyset\wedge e_{\{2\}}\wedge e_{\{1,3\}}\wedge e_{\{1,2,3\}}\!-
e_\emptyset\wedge e_{\{1\}}\wedge e_{\{2,3\}}\wedge e_{\{1,2,3\}}\!+
$$
$$
e_{\{2\}}\wedge e_{\{3\}}\wedge e_{\{1,2\}}\wedge e_{\{1,3\}}\!-
e_{\{1\}}\wedge e_{\{3\}}\wedge e_{\{1,2\}}\wedge e_{\{1,3\}}\!+
e_{\{1\}}\wedge e_{\{2\}}\wedge e_{\{1,3\}}\wedge e_{\{2,3\}}\!+
2e_{\{1\}}\wedge e_{\{2\}}\wedge e_{\{3\}}\wedge e_{\{1,2,3\}}\!.
$$  
Hence this is the $SO(7)$-invariant in $\wedge^4\Gamma$. 

The $SO(6)$-subrepresentations in $S$ 
(where $SO(6)$ acts on $\langle e_1,\ldots,e_6\rangle \subset\CC^7$) are
$$
W\cong \wedge^{even} R,\qquad W^*\cong  \wedge^{odd} R.
$$
The $so(7)$-invariant given above has a non-trivial component
$e_\emptyset\wedge e_{\{1,2\}}\wedge e_{\{1,3\}}\wedge e_{\{2,3\}}\in \wedge^4W$, hence the theorem follows.
\qed

\


\

\end{document}